\magnification = 1100

\input amstex

\documentstyle{amsppt}

\topmatter
\title Smooth norms and  approximation in Banach spaces of
the type $\Cal C(K)$
\endtitle

\author Petr H\'ajek and Richard Haydon \endauthor
\address Mathematical Institute, Czech Academy of Science,
\v Zitn\'a 25,
Praha, 11567, Czech rep
\endaddress
\email  hajek\@math.cas.cz\endemail
\address Brasenose College, Oxford, OX1 4AJ, England
\endaddress
\email richard.haydon\@brasenose.oxford.ac.uk \endemail
\subjclass 46B03, 46B26 \endsubjclass
\thanks{The research presented in this article was supported by
grants AV 1019205  and GA CR 201/04/0090 awarded by the Czech
Academy of Sciences.}\endthanks

\endtopmatter
\document
\heading
1. Introduction
\endheading

We prove two results about smoothness in Banach spaces of the type
$\Cal C (K)$. Both build upon earlier papers of the first named author [3,4].
We first establish a special case of a conjecture that
remains open for general Banach spaces and concerns smooth
approximation.  We recall that a {\it bump function} on a Banach
space $X$ is a function $\beta:X\to \Bbb R$ which is not
identically zero, but which vanishes outside  some bounded set.
The existence of bump function of class $\Cal C^1$  implies that
the Banach space $X$ is an {\it Asplund space}, which, in the case
where $X=\Cal C(K)$, is the same as saying that $K$ is {\it
scattered}.  It is a major unsolved problem to determine whether
every Asplund space has a $\Cal C^1$ bump function.  Another open
problem is whether the existence of just one bump function of some
class $\Cal C^m$ on a Banach space $X$
implies that all continuous functions on $X$ may be uniformly
approximated by functions of class $\Cal C^m$. It is to this
question that we give a positive answer (Theorem~2) in the special
case of $X=\Cal C(K)$.

Our second result represents some mild progress with a conjecture
made by the second author in [7].  The analysis in that paper of compact spaces
constructed using trees suggested that for a compact space $K$ the
existence of an equivalent {\it norm} on $\Cal
C(K)$ which is of class $\Cal C^1$ (except at 0 of course) might
imply the existence of such a norm which is of class
$\Cal C^\infty$. Certainly, this is what happens with norms
constructed using linear {\it Talagrand operators} as in [5,6,7].
The other important (and older) method of obtaining $\Cal C^1$
norms is to construct a norm with locally uniformly rotund dual
norm.  What we show in Theorem~2 is that, whenever $\Cal C(K)$ admits
an equivalent norm with LUR dual norm, then there is also an
infinitely differentiable equivalent norm on $\Cal C(K)$.

For background on smoothness and renormings in Banach spaces,
including an account of Asplund spaces, we
refer the reader to [1].  In particular, an account is given there
of the connection between smooth approximability of continuous
functions and the existence of smooth {\it partitions of unity}.
Following what seems to be standard practice in the literature, we
have chosen to state the formal version of our first theorem
(Theorem 1) in terms of partitions of unity, rather than
approximation.

Formalizing a definition which appears implicitly in [6,7], we shall
say that a mapping $T:\Cal C(K)\to c_0(K\times \Gamma)$ is a
(non-linear) Talagrand operator of class $\Cal C^m$ if
\roster
\item for each non-zero $f \in \Cal C(K)$ there exist $t\in K,
u\in \Gamma$ such that $|f(t)|=\|f\|_\infty$ and $(Tf)(t,u)\ne 0$;
\item each coordinate function $f\mapsto (Tf)(t,u)$ is of class
$\Cal C^m$ on the set where it is not zero.
\endroster

It follows from Corollary~3 of [6] that if $\Cal C(K)$ admits a
Talagrand operator of class $\Cal C^m$ then $\Cal C(K)$ has a
bump function of the same class.  Corollary~2 of [6], or Theorem 2
of [5], shows that if $\Cal C(K)$ admits a {\it linear} Talagrand
operator then $\Cal C(K)$ admits an equivalent $\Cal C^\infty$
norm.
Whilst there certainly exist examples of compact $K$ such that
$\Cal C(K)$ has a $\Cal C^\infty$ renorming but no linear
Talagrand operator (for instance, the Cieselski-Pol space, [1,
p260]), it is worth noting that, by the first of the theorems in
this paper, a non-linear Talagrand operator exists whenever there
is a bump function.

\heading
Admissible subsets of compact scattered spaces
\endheading

In this short section, we establish some notation that will be
used consistently later on, and prove one easy result.
We consider a compact scattered space $K$.  The derived set $K'$
is defined as usual to be the set of points $t$ in $K$ which are not
isolated in $K$. Successive derived sets $K^{\alpha)}$ are defined
by the transfinite recursion
$$
K^0=K; \quad K^\beta= \bigcap_{\alpha<\beta} (K^{\alpha})'.
$$
There is an ordinal
$\delta$ such that $K^{(\delta)}$ is non-empty and finite (so that
$K^{(\delta+1)}=\emptyset$).  For each $t\in K$ there is a unique
ordinal $\alpha(t)\le \delta$ such that $t\in
K^{\alpha(t)}\setminus K^{\alpha(t)+1}$.  Since $t$ is an isolated
point of $K^{\alpha(t)}$, there is a compact open subset $V$ of
$K$ such that $V\cap K^{\alpha(t)}=\{t\}$; we choose such a $V$
and call it $V_t$.  For finite subsets $B$ of $K$, we set $V_B=
\bigcup_{t\in B}V_t$.

\proclaim{Lemma 1} Let $B$ be a nonempty finite subset of $K$ and let
$\alpha=\alpha(B)$ be maximal subject to $B\cap K^{(\alpha)}\ne \emptyset\}$. Then
$V_B \cap K^{\alpha}=B\cap K^{\alpha}$ and hence $V_B\cap
K^{\alpha+1}=\emptyset$.
\endproclaim
\demo{Proof}
Let $t$ be in $B$.  If $\alpha(t)<\alpha$ then $V_t\cap
K^{(\alpha)}=\emptyset$, whilst if $\alpha(t)=\alpha$ then $V_t\cap
K^{(\alpha)}=\{t\}$.  Thus $V_B \cap K^{\alpha}=B\cap K^{\alpha}$,
as claimed.
\enddemo

We shall say that a finite subset $A$ of $K$ is {\it admissible}
if $s\notin V(t)$ whenever $s$ and $t$ are distinct elements of
$A$.

\proclaim{Lemma 2} Let $K$ be a compact scattered space and let $H$
be a nonempty, closed subset of $K$.  There is a unique admissible
set $A$ with the property that $A\subseteq H\subseteq V(A)$.
\endproclaim
\demo{Proof}
We start by describing a recursive procedure which constructs one
possible admissible $A$ with the required property.  Let
$\alpha_0=\max\{\alpha: H\cap K^{(\alpha)}\ne \emptyset\}$; thus,
$H\cap K^{(\alpha_0)}$ is a nonempty finite set, which we shall
call $A_0$.  If $H\subseteq V_{A_0}$ we set $A=A_0$ and stop.
Otherwise, we set $H_1=H\setminus V_{A_0}$, $\alpha_1 = \max
\{\alpha: H_1\cap K^{(\alpha)}\ne \emptyset\}$, $A_1=H_1\cap
K^{(\alpha_1)}$, and continue.  In this way, we construct a
decreasing (and so, necessarily finite) sequence
$\alpha_0>\alpha_1>\cdots>\alpha_l$ of ordinals, and finite sets
$A_j= H\cap K^{(\alpha_j)}\setminus V_{A_0\cup \dots \cup
A_{j-1}}$, in such a way that $H\subseteq V_{A_0\cup \dots \cup
A_l}$.  By construction, the set $A=A_0\cup \dots \cup
A_l$ is admissible.

We now show uniqueness. It will be convenient to proceed by transfinite
induction on $\alpha_0$.  Let $B$ be admissible and suppose
that $B\subseteq H\subseteq V_B$. By Lemma~1, $\alpha(B)=
\alpha_0$ and $B\cap K^{(\alpha_0)} \subseteq C\cap K^{(\alpha_0)} \subseteq
V_B\cap K^{(\alpha_0)} = B\cap K^{(\alpha_0)}$.  Thus
$A_0\subseteq B$.  We now have a closed set $H_1=H\setminus V_{A_0}$ and an
admissible set $B_1=B\setminus A_0$ with $B_1\subseteq H_1\subseteq
V_{B_1}$.  Since $\alpha_1 =\max\{\alpha: H_1\cap K^{(\alpha)}\ne
\emptyset\}<\alpha_0$, we may use our inductive hypothesis to
deduce that $B_1 = A\setminus A_0$, whence $B=A$.
\enddemo

\heading
Improving bump functions on $\Cal C(K)$
\endheading

Let $X$ be a Banach space that admits a bump function of class
$\Cal C^m$; so there is a function $\alpha\in \Cal C^m(X)$ such
that $\alpha(0)=1$ while $\alpha(x)=0$ for $\|x\|\ge 1$.  By forming
$\beta$, where $\beta(x) = \phi(\alpha(x/R))$, with $R>0$ and $\phi:\Bbb R\to
\Bbb R$ suitably chosen, we obtain a function of class $\Cal C^m$,
taking values in $[0,1]$
satisfying
$$
\align
\beta(x) &=1 \qquad\text{when $\|x\|\le 1$}\\
\beta(x) &=0 \qquad\text{when $\|x\|\ge R$}.
\endalign
$$
Of course, if $X$ admits partitions of unity of class $\Cal C^m$,
then (starting with a partition of unity each of whose members has
support of diameter at most $\epsilon$) we easily obtain a
function $\beta$ satisfying the above conditions, with
$R=1+\epsilon$ and $\epsilon$ an arbitrarily small positive real
number. In general,  we do not know whether a bump function can
always be ``improved'' in this way.  This section is devoted to
showing how to achieve such an improvement in the case where $X$
is a space $\Cal C(K)$ equipped with the supremum norm.  We start
with an elementary and no doubt well known exercise in calculus.

\proclaim{Lemma 3}
Let $K$ be a compact space and let $\theta:\Bbb R \to \Bbb R$ be
of class $\Cal C^m$.  Then the mapping $\Theta:\Cal C(K)\to \Cal
C(K)$ given by $\Theta(f)=\theta\circ f$ is of class $\Cal C^m$.
\endproclaim
\demo{Proof}
We proceed by induction on $m$.  For $m=0$ we are merely assuming
$\theta$ to be continuous, and continuity of $\Theta$ follows from
the uniform continuity of $\theta$ on bounded subsets of $\Bbb R$.

If $m\ge 1$ we consider $f,h$ in $\Cal C(K)$ and apply the mean
value theorem point by point to obtain
$$
\theta(f(t)+h(t))-\theta(f(t)) = \theta'(f(t) + \zeta (t) h(t))h(t),
$$
with $0<\zeta(t)<1$.  The uniform continuity of $\theta'$ on
bounded subsets of $\Bbb R$ now tells us that the right hand side
of the above equality equals $\theta'(f(t))+\text
o(\|h\|_\infty)$. So $\Theta$ is differentiable with
$$
\text D\Theta(f)\cdot h = (\theta'\circ f)\times h.
$$
The linear mapping $\text D\Theta(f): \Cal C(K)\to \Cal C(K)$ is thus
the operator $\text M_{\theta'}$ of multiplication by $\theta'\circ f$.  Thus the
derivative $\text D \Theta:\Cal C(K)\to \Cal L(\Cal C(K))$ may be factored as follows
$$
\Cal C(K)\to \Cal C(K) \to \Cal L(\Cal C(K)),
$$
where the first factor is $f\mapsto \theta'\circ f$ and the second
is the linear isometry $g\mapsto \text M_g$.  Our inductive
hypothesis tells us that the first of these is of class $\Cal
C^{m-1}$.  So $D\Theta$ is of class $C^{m-1}$ and $\Theta$ of class $C^{m}$.
\enddemo

\proclaim{Proposition 1}
Let $K$ be a compact space such that $\Cal C(K)$ admits a bump
function of class $\Cal C^m$.  Then, for all real numbers
$\eta>\xi>0$, there is a function $\beta_{\xi,\eta}:\Cal C(K)\to
[0,1]$ of class $\Cal C^m$ such that
$$
\beta_{\xi,\eta}(f)= \cases 1\quad\text{when $\|f\|_\infty\le
\xi$}\\
0\quad\text{when $\|f\|_\infty\ge
\eta$}
\endcases
$$
\endproclaim
\demo{Proof}
By hypothesis, there exists a function $\alpha:\Cal C(K)\to \Bbb
R$, of class $\Cal C^m$, such that $\alpha(0)=1$ while $\alpha
(f)=0$ for $\|f\|_\infty\ge 1$.  As in our introductory remarks,
we may assume that $\alpha$ takes values in $[0,1]$.
We define $\beta_{\xi,\eta}$ by
$$
\beta_{\xi,\eta}(f) = \alpha(\theta\circ f),
$$
where $\theta:\Bbb R\to \Bbb R$ is a function of class $\Cal
C^\infty$ chosen so that
$$
\theta(x) = \cases 0\quad\text{when $|x|\le \xi$}\\
1\quad\text{when $|x|\ge \eta$}
\endcases
$$
\enddemo

\heading
Construction of a Talagrand operator and of partitions of unity
\endheading

This section is devoted to a proof of the following theorem.

\proclaim{Theorem 1}
Let $K$ be a compact space and let $m$ be a positive integer or $\infty$.
The following are equivalent:
\roster
\item $\Cal C(K)$ admits a bump function of class $\Cal C^m$;
\item $\Cal C(K)$ admits a Talagrand operator of class $\Cal C^m$;
\item $\Cal C(K)$ admits a partitions of unity of class $\Cal C^m$.
\endroster
\endproclaim

It will be enough to prove that (1) implies both (2) and (3).  We
start by showing how to construct a Talagrand operator, starting
with a bump function on $\Cal C(K)$.  As we remarked in the Introduction, the
existence of a smooth bump function forces $K$ to be scattered.  So we can
use the notion of admissible set as developed above.  We let $\Cal Q$ be the set
of all triples $(\xi, \eta, \zeta)$ in $\Bbb Q^3$ with $0<\xi<\eta<\zeta$, and
write $\Cal A$ for the set of all admissible subsets of $K$.
Choose positive real numbers $c(\xi, \eta,\zeta)$ with $\sum_{(\xi,
\eta,\zeta)\in \Cal Q} c(\xi, \eta,\zeta)<\infty$.  For
$0<\xi<\eta$ let $\beta_{\xi,\eta}$ be as in Proposition~1, and,
finally, for $0<\eta<\zeta$, let $\phi_{\eta,\zeta}:\Bbb R\to [0,1]$ be of
class $\Cal C^\infty$ with
$$
\phi_{\eta,\zeta} (x) = \cases 0\quad\text{when $|x|\le \eta$}\\
1\quad\text{when $|x|\ge \zeta$}
\endcases
$$

We define $T:\Cal C(K)\to \ell^\infty(K\times \Cal Q\times \Cal A)$
by
$$
(Tf)(s,\xi,\eta,\zeta,A) = c(\xi,\eta,\zeta)
\beta_{\xi,\eta}(f\times \chi_{K\setminus V_A})\prod_{t\in
A}\phi_{\eta,\zeta}(f(t))\chi_A(s).
$$
We notice that, for this expression to be non-zero, we need
$A\subseteq F\subseteq V_A$, where $F$ is the closed set $\{t\in
K: |f(t)|\ge \eta\}$.  Now, we know by Lemma~2 that, for given $\eta$ and $f$, there is
just one $A$ for this is true.  It follows easily that $T$ takes
values in $c_0(K\times \Cal Q\times \Cal A)$.  It is also clear
that each coordinate of $T$, that is to say, each mapping
$f\mapsto (Tf)(s,\xi,\eta,\zeta,A)$, is of class $\Cal C^m$.

To show that $T$ has the Talagrand property, we consider $f\ne 0$
and set $F=\{t\in K:|f(t)|=\|f\|_\infty\}$.  Let $A$ be the
admissible set for which $A\subseteq F\subseteq V_A$ and choose
rationals $0<\xi<\eta<\zeta$ such that $\zeta<\|f\|_\infty$ and
$\xi>\|f\times \chi_{K\setminus V_A}\|_\infty$.  For any $s\in A$
we have $|f(s)|=\|f\|_\infty$ and $(Tf)(s,\xi,\eta,\zeta)\ne 0$.

We now pass to the construction of partitions of unity.  We shall
proceed by transfinite recursion on the derived length of $K$.
Recall from Section~2 that there is an ordinal $\delta$ such that
$K^{(\delta)}$ is finite and non-empty, so that $K^{(\delta+1)}$
is the first empty derived set of $K$.  We assume inductively
that,
if $V$ is a compact space with $V^{(\delta)}=\emptyset$ and such
that $\Cal C(V)$ has a bump function of class $\Cal C^m$, then
$\Cal C(V)$ admits $\Cal C^m$ partitions of unity.  We need to
show that $\Cal C(K)$ also admits $\Cal C^m$ partitions of unity.
To do this it will be enough to construct partitions of unity on
the finite-codimensional subspace $X=\{f\in \Cal C(K): f(t)=0 \text{ for all } t\in
K^{(\delta)}$.  We shall use the following result from an earlier
paper by the second author.

\proclaim{Proposition 2 {\smc(Theorem 2 of [6])}}
Let $X$ be a Banach space, let $L $ be a set and let $m$ be a positive
integer or $\infty$. Let $T:X\to c_0(L
)$ be a function such that each coordinate $x\mapsto T(x)_\gamma $ is of class
$\Cal C^m$ on the set where it is non-zero. For each finite subset $F$ of $L $,
let $R_F:X\to X$ be of class $\Cal C^m$ and  assume that the following hold:
\roster
\item
for each $F$, the image $R_F[X]$ admits $\Cal C^k$ partitions of unity;
\item
$X$ admits a $\Cal C^k$ bump function;
\item
for each $x\in X$ and each $\epsilon >0$ there exists $\lambda >0$ such that
$\|{x-R_Fx}\|<\epsilon $ if we set $F =\{u \in L
:|(Tx)(u )|\ge \lambda \}$.
\endroster
Then $X$ admits $\Cal C^m$ partitions of unity.
\endproclaim

In applying this result we shall take $L$ to be $K\times \Cal
Q\times \Cal A_0$, where $\Cal A_0$ consists of the admissible
subsets $A$ such that $A\cap K^{(\delta)}=\emptyset.$ The operator
$T$ is the Talagrand operator constructed above (though with the
argument $A$ restricted to lie in $\Cal A_0$). We have already shown that
$T$ takes values in $c_0(L)$ and that the coordinates of $T$ are of class $\Cal C^m$.
We define the
reconstruction operators $R_F$ as follows:  if $F\subset L$ has
elements $(s_i, \xi_i, \eta_i,\zeta_i, A_i)$  $(0\le i<n)$, we set
$V(F)=\bigcup_{i<n} V_{A_i}$ and define $R_F(f) = f \times
\chi_{V(F)}$.  So $R_F:X\to X$ is a bounded linear operator and
the image $R_F$ may be identified with $\Cal C(V(F))$, which, by
our inductive hypothesis, admits partitions of unity of class
$\Cal C^m$.

It only remains to check that (3) holds, so let $f\in \Cal C(K)$
and $\epsilon>0$ be given.  Let $H$ be the set $\{t\in K:
|f(t)|\ge \epsilon\}$ and let $A$ be the admissible set such that
$A\subseteq H\subseteq V_A$.  For suitably chosen
$0<\xi<\eta<\zeta<\epsilon$ we have
$$
(Tf)(s,\xi,\eta,\zeta,A) = c(\xi,\eta,\zeta)>0
$$
for all $s\in A$.  We set $\lambda= c(\xi,\eta,\zeta)$ and note that
$V(F)\supseteq V_A$.  So
$$
\align
\|f-R_Ff\|_\infty &= \|f\times \chi_{K\setminus V(F)}\|_\infty\\
    &\le \|f\times \chi_{K\setminus V_A}\|_\infty <\epsilon.
    \endalign
$$

\heading
Infinitely differentiable norms
\endheading

We shall now prove the second theorem of this paper.

\proclaim{Theorem 2}
Let $K$ be a compact space such that $\Cal C(K)$ admits an
equivalent norm with locally uniformly rotund dual norm.  Then $\Cal
C(K)$ admits an equivalent which is of class $\Cal C^\infty$ on
$X\setminus \{0\}$.
\endproclaim

The norm which we construct will be a {\it generalized Orlicz
norm}, defined on the whole of $\ell^\infty(K)$, which we shall
show to be infinitely smooth on the subspace $\Cal C(K)$. We
recall some definitions.  Suppose that, for each $t\in K$, we are
given a convex function $\phi_t= \phi(t,\cdot):[0,\infty)\to [0,\infty)$
satisfying $\phi(t,0)=0$, $\lim_{x\to\infty} \phi(t,x)=\infty$
(that is to say an {\it Orlicz function}).  The generalized Orlicz
space $\ell_{\phi(\cdot)}(K)$ is defined to be the space of all
functions $f:K\to \Bbb R$ such that $\sum_{t\in
K}\phi(t,|f(t)|/\rho)<\infty$ for some $\rho\in (0,\infty)$.  The
generalized Orlicz norm of such a function is defined to be
$$
\|f\|_{\phi(\cdot)} = \inf\{\rho>0: \sum_{t\in
K}\phi(t,|f(t)|/\rho)\le 1\}.
$$

The first of the following lemmas is elementary and the second
uses the familiar idea of ``local dependence on finitely many
coordinates.''

\proclaim {Lemma 4}
Suppose that there exist positive real numbers $R<S$ such that
$\phi(t,R)=0$ and $\phi(t,S)\ge 1$ for all $t\in K$.  Then
$\ell_{\phi(\cdot)}(K)=\ell^\infty(K)$ and
$$
R\|f\|_{\phi(\cdot)}\le \|f\|_\infty \le S\|f\|_{\phi(\cdot)}.
$$
\endproclaim

\proclaim{Lemma 5}
Let $\phi(\cdot)$, $R$ and $S$ be as in the preceding lemma, and
let $X$ be a linear subspace of $\ell^\infty(K)$.  Suppose that,
whenever $f\in X$ and $\|f\|_{\phi(\cdot)}=1$, there exists a
positive real number $\delta$ and a finite subset $F$ of $K$ such that
$g(t)=0$ whenever $g\in X$, $\|f-g\|_\infty<\delta$ and $t\notin
F$.  Assume further that each of the functions $\phi(t,cdot)$ is
of class $\Cal C^\infty$.  Then the generalized Orlicz norm
$\|\cdot\|_{\phi(\cdot)}$ is of class $\Cal C^\infty$ on
$X\setminus \{0\}$.
\endproclaim
\demo{Proof}
If $f$, $F$ and $\delta$ are as in the statement of the lemma,
then the function $\Phi$ defined by
$$
\Phi(g) = \sum_{t\in K}\phi(t,|g(t)|)
$$
is of class $\Cal C^\infty$ on $\{g\in X:\|f-g\|_\infty
<\delta\}$, since it coincides with the finite sum $\sum_{t\in
F}\phi(t, |g(t)|)$ of given $\Cal C^\infty$ functions.  Thus our
hypothesis tells that there is an open subset $U$ of $X$
containing the set $\{f\in X: \|f\|_{\phi(\cdot)}=1\}$ and such
that $\Phi$ is $\Cal C^\infty$ on $U$.  We define $V=\{h,\rho)\in
(X\setminus \{0\})\times (0,\infty): \rho^{-1}h\in U$, an open set
in $X$.  On $V$ we define $\Psi(h,\rho)= \Phi (\rho^{-1}h)$,
which is of class $\Cal C^\infty$.  For each $h\in X\setminus
\{0\}$ there is a unique $\rho= \|h\|_{\phi(\cdot)}$ such that
$(h,\rho)\in V$ and $\Psi(h,\rho)=1$.  Moreover, we may calculate
the partial derivative
$$
\text D_2\Psi(h, \rho) = -\rho^{-2}\sum_t \phi'_t( \rho^{-1}|h(t)|),
$$
and note that this is nonzero when $\rho = \|h\|_{\phi(\cdot)}$,
since $\phi'_t(x)>0$ whenever $\phi_t(x)>0$.  Thus the Implicit
Function Theorem may be applied to conclude that
$\|\cdot\|_{\phi(\cdot)}$ is of class $\Cal C^\infty$ on $X\setminus
\{0\}$.
\enddemo

To choose suitable Orlicz functions $\phi_t$ in our theorem, we
shall need to use the special properties of $K$.
The assumption that $\Cal C(K)^*$ has an equivalent LUR dual norm
implies (and, by a theorem of Raja [9], is actually equivalent to) the
compact space $K$ being $\sigma$-discrete.  So we may assume that
there are pairwise disjoint subsets $D_i$ of $K$, each one discrete in its subspace
topology, with $K=\bigcup_{i\in \omega} D_i$. We note in passing
that we are not assuming the $D_i$ to be closed, merely that $D_i$
has empty intersection with its derived set $D_i'$.  We fix positive
real numbers $r_i<1$ with $\prod _{i\in \omega}r_i>0$ and, for
$t\in K$, define two real numbers
$$
\align
\alpha(t) &= \prod\{ r_i: t\in \bar D_i\}\\
\beta(t) &= \prod\{ r_i: t\in D_i'\}.
\endalign
$$
We notice that $\beta(t) = \alpha(t)\times r_j$ where $j$ is the
(unique) natural number such that $t\in D_j$. (Note that it is here
(and only here) that we use the discreteness hypothesis that $D_j\cap
D_j'=\emptyset$.)
In particular, therefore,
$0<\alpha(t)<\beta(t)<1$.  So we may choose an infinitely differentiable
Orlicz function $\phi_t$ such that
$$
\align
\phi_t(x) &= 0 \quad \text{when } x\le \alpha(t)\\
\phi_t(x) &> 1 \quad \text{when } x\ge \beta(t).
\endalign
$$

We are going to show that Lemma~5 may be applied to these Orlicz
functions and the subspace $\Cal C(K)$ of $\ell^\infty$.
It is convenient to state one of the ingredients of this proof as
a property of the functions $\alpha$ and $\beta$.

\proclaim{Lemma 6}
Let $t_n$ be a sequence of distinct elements of $K$ which converges to some
$t\in K$.  Then $\beta(t)\le \liminf \alpha(t_n)$.
\endproclaim
\demo{Proof}
By taking subsequences and diagonalizing, we may assume that
$\alpha(t_n)$ tends to a limit as $n\to \infty$, and also that,
for each $i\in \omega$, either all of $t_{i+1},
t_{i+2},\dots$ are in $\bar D_i$, or else none is.  Let $M$ be the
set of natural numbers $i$ such that $t_{i+1},\dots$ are in $\bar
D_i$.  Then, for each $n$ we have
$$
\prod_{i<n,\ i\in M} r_i\times \prod _{i\ge n} r_i \le
\alpha(t_n) \le \prod _{i\in M} r_i,
$$
whence $\alpha(t_n)\to \prod_{i\in M}r_i$ as $n\to \infty$.
On the other hand, since the $t_n$ are distinct and
$t_n\in D_i$ whenever $n>i\in M$, it must be that the limit point
$t$ is in the derived set $D'_i$ whenever $i\in M$.  Thus
$$
\beta(t) = \prod\{ r_i: t\in D_i'\}\le \prod_{i\in M}r_i.
$$
\enddemo

To complete the proof of the theorem, we consider $f\in \Cal C(K)$
with $\|f\|_{\phi(\cdot)}=1$.  If no $\delta $ and $F$ exist with
the property of Lemma~5, there exist a sequence $(f_n)$ in $\Cal
C(K)$ converging uniformly to $f$ and a sequence of distinct
elements $(t_n)$ of $K$ such that $\phi(t_n,|f_n(t_n)|)>0$ for all
$n$.  For this to be the case, if must be that $|f_n(t_n)|\ge
\alpha(t_n)$. Extracting a subsequence, we may suppose that the
sequence $(t_n)$ converges to some $t\in K$.  Now by uniform
convergence and the continuity of
$f$ we have $f(t) = \lim_n f_n(t_n)$, so that $|f(t)|\ge
\beta(t)$ by Lemma~6.  Thus $\phi(t,|f(t)|)>1$ and
$\|f\|_{\phi(\cdot)}>1$, a contradiction, which ends the proof.

\heading
Remarks
\endheading
It follows from Lemma~4 that the norm constructed in Theorem~2
satisfies
$$
\|f\|_\infty \le \|f\|_{\phi(\cdot)}\le \alpha^{-1} \|f\|_\infty,
$$
where $\alpha= \prod_{i\in \omega} r_i$.  Since we may arrange for $\alpha$ to
be arbitrarily close to 1,  we have shown
that the supremum norm may be uniformly approximated on bounded
subsets of $\Cal C(K)$ by infinitely differentiable norms.  We do
not know whether all equivalent norms on $\Cal C(K)$ (with $K$
$\sigma$-discrete) can be thus approximated.  In particular, we do
not know this for the locally uniformly rotund norm recently
constructed by the second author [8].
Of course it follows from the results of the present paper that any equivalent norm on
our $\Cal C(K)$ (like any other continuous function)  can be
uniformly approximated on bounded subsets by infinitely
differentiable {\it functions}.

\enddocument